\documentclass[preprint,12pt]{elsarticle}

\usepackage[fleqn]{amsmath}
\usepackage[]{amsfonts}
\usepackage[]{amsthm}
\usepackage[]{amssymb}
\usepackage{empheq}
\usepackage{eucal}

\usepackage{graphics}
\usepackage{graphicx}
\usepackage{commath}
\usepackage{mathtools}
\usepackage{tabu}
\newcolumntype{K}[1]{>{\centering\arraybackslash}m{#1}}

\usepackage[utf8]{inputenc}

\usepackage{setspace}
\usepackage{titlesec}
\titlespacing{\section}{0pt}{15pt}{5pt}
\titlespacing{\subsection}{0pt}{10pt}{3pt}
\titlespacing{\subsubsection}{0pt}{5pt}{1pt}

\usepackage{bm}

\usepackage{array}
\usepackage{pdflscape}
\usepackage{caption}

\newcommand{\stkout}[1]{\ifmmode\text{\sout{\ensuremath{#1}}}\else\sout{#1}\fi}

\usepackage[table]{xcolor}

\expandafter\def\expandafter\normalsize\expandafter{%
    \normalsize
    \setlength\abovedisplayskip{5pt}
    \setlength\belowdisplayskip{5pt}
    \setlength\abovedisplayshortskip{10pt}
    \setlength\belowdisplayshortskip{10pt}
}

\usepackage{lineno}

\journal{Transportation Research Part B}

\begin{document}

\begin{frontmatter}




\title{Dynamic Traffic Assignment using the Macroscopic Fundamental Diagram: A Review of Vehicular and Pedestrian Flow Models}

\author[GATECH]{Rafegh Aghamohammadi}
\author[GATECH]{Jorge A. Laval\corref{cor}}  
\ead{jorge.laval@ce.gatech.edu}

\cortext[cor]{Corresponding author} 

\address[GATECH]{School of Civil and Environmental Engineering, Georgia Institute of Technology}

\begin{abstract}

Traditional DTA models 
of large cities suffer from prohibitive computation times and calibration/validation can become major challenges faced by practitioners. 
The empirical evidence in 2008 in support of the existence of a Macroscopic Fundamental Diagram (MFD) on urban networks 
led to the formulation of discrete-space models, where the city is divided into a collection of reservoirs.
Prior to 2008, a large body of DTA models based on pedestrian flow models had been formulated in continuum space as 2-dimensional conservation laws where the speed-density relationship can now be interpreted  as the MFD.  
Perhaps surprisingly, we found that this continuum-space literature has been mostly unaware of MFD theory, and  no attempts exist to verify  the assumptions of MFD theory. This has the potential to create significant inconsistencies, and research is needed to analyze their extent and ways to resolve them.   We also find that further research is needed to 
(i) incorporate departure time choice,
(ii) improve existing numerical methods,  possibly extending  recent advances on the one-dimensional kinematic wave (LWR) model,
(iii) study the properties of system optimum solutions,
(iv) examine the real-time applicability of current continuum-space models compared to traditional DTA methods, and
(v) formulate anisotropic models for the interaction of intersecting flows.

\end{abstract}


\end{frontmatter}

\section{Introduction}

Dynamic traffic assignment (DTA) models based on a link-level representation of the network are becoming
the mainstream method for urban planning, either under user equilibrium (UE) or system optimum (SO) objectives. Despite the vast body of literature in this field, there are still some important drawbacks. Arguably, the main challenge faced by practitioners is the prohibitive computation times necessary to achieve equilibrium on large-scale networks, which typically have hundreds of thousands of links and nodes. This makes calibration and validation of such large networks a daunting task.

\begin{figure} [!b]
\begin{center}
\includegraphics[width=13cm]{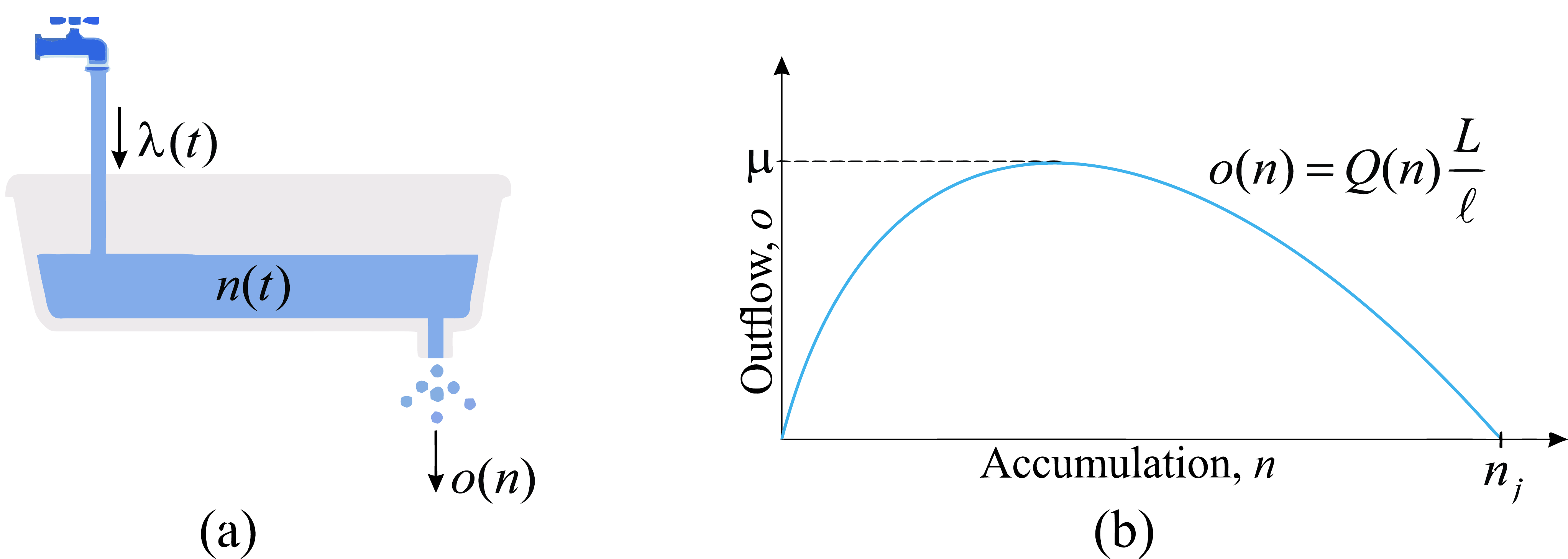}
\vspace{-10pt}
    \caption{Illustration of: (a) reservoir model and (b) resulting o-MFD diagram.}
    \vspace{-10pt}
    \label{f0}
\end{center}
\end{figure}

In the past 50 years, various theories have been proposed to describe the traffic flow of urban networks at an aggregate level \citep{Smeed1967Road,Herman1979Two,Herman1984Characterizing}. The recent empirical verification of the existence of a network-level Macroscopic Fundamental Diagram (MFD) on congested urban areas in 2008 opened up a new paradigm \citep{Daganzo2007Urban,Geroliminis2008Existence}. The MFD
\begin{equation}\label{eqMFD}
Q=Q(n), \hskip 4.9cm \mbox{(average flow MFD)}
\end{equation}
\noindent
gives the average flow $Q$ on network as a function of the number of vehicles inside the network, $n$, arguably independently of trip origins and destinations, and route choice. 
This makes the MFD an invaluable tool to overcome the difficulties of traditional planning models. For example, on a single region one can construct a reservoir type model to approximate the evolution of the accumulation of vehicles inside the reservoir. We can then estimate any average traffic variable of interest, provided that we know the outflow MFD, o-MFD, defined as
\begin{equation}\label{eqOMFD}
o(n)=Q(n)\frac{L}{\ell}, \hskip 5.7cm \mbox{(o-MFD)}
\end{equation}
where $L$ corresponds to the length of network and  $\ell$ is the trip length, typically assumed identical for all commuters. The o-MFD, $o(n)$, gives the number of trip completions per unit time, as a function of $n$. The simplest model for the traffic dynamics inside the network is the simple reservoir (or bathtub) model, as shown in Fig. \ref{f0}(a), which simply states the conservation of vehicles via the following ordinary differential equation (ODE):
\begin{subequations}\label{model0}
\begin{empheq}[left=\textbf{}\empheqlbrace]{align}
	n '(t)&=\lambda (t)-o (n), & \mbox{(reservoir dynamics)}\label{mfdode}\\
	n (0)&=n_0, & \mbox{(initial conditions)}\label{boundary}
\end{empheq}
\end{subequations}
\noindent
where $\lambda(t)$ is the demand inflow into the network at time $t$,  primes denote differentiation and  $n_0$ specify the initial conditions. Fig. \ref{f0}(b) illustrates a typical o-MFD.

The shape of MFD depends on network topology and control parameters such as block length, existence of turn-only lanes, and traffic light settings. 
Up until recently,
estimating the MFD was not an easy task because the method of cuts in \citep{Daganzo2008Analytical}
becomes intractable for real-life networks and one needs to resort to simulation methods, which defeats the purpose of macroscopic modeling.
However, \citet{Laval2015Stochastic} show that (the probability distribution of) the MFD can be well approximated by a function of mainly two  parameters: the density of traffic lights and the mean red to green ratio across the network. Since these are observable parameters, we now have a simple method for estimating the MFD on arbitrary road networks.

We conclude that it is now possible to formulate  DTA models based on a macroscopic representation of the network. This representation may be expressed in (i) continuum  space, where at each $(x,y)$-point there are macroscopic functions that give the MFD, the demand, the density of traffic lights, etc.; in (ii) discrete space,
where the modeling region is divided in a finite number of zones describing a tessellation of the $(x,y)$-plane. Each zone has a well-defined MFD with traffic dynamics given by the conservation ODE,   \eqref{model0}. 

It turns out that there is a large body of related literature prior to 2008, when the MFD was verified empirically, that dealt with similar problems. For example, we have DTA models in continuum space that consist of a conservation law in two dimensions, supplemented with a Hamilton-Jacobi equation for the minimum paths. Mathematically, these models are identical to the continuum pedestrian flow models in the literature, which is why we have included pedestrian literature here. All of these models rely on the existence of a speed-density equilibrium relationship, which could be thought of as the MFD. 

Our purpose here is to review the literature in light of the aforementioned considerations. In particular, we would like to identify the current challenges to bringing to practice efficient macroscopic DTA models for cities. Towards this end, the remainder of this paper is organized as follows. 
Section 2 describes the formulation of continuum-space models, and section 3 presents their solution methods, including analytical and numerical. Section 4 discusses discrete-space models, and finally conclusions and outlook for further research are presented in section 5.

It is worth to mention that in order to keep the consistency of the notation through this paper, we have changed the notation of some of the models presented here. Also, a \textbf{bold-typed reference} in this paper indicates the beginning of its review.

\section{Continuum-space models}

\begin{figure} [h]
\begin{center}
\includegraphics[width=110mm]{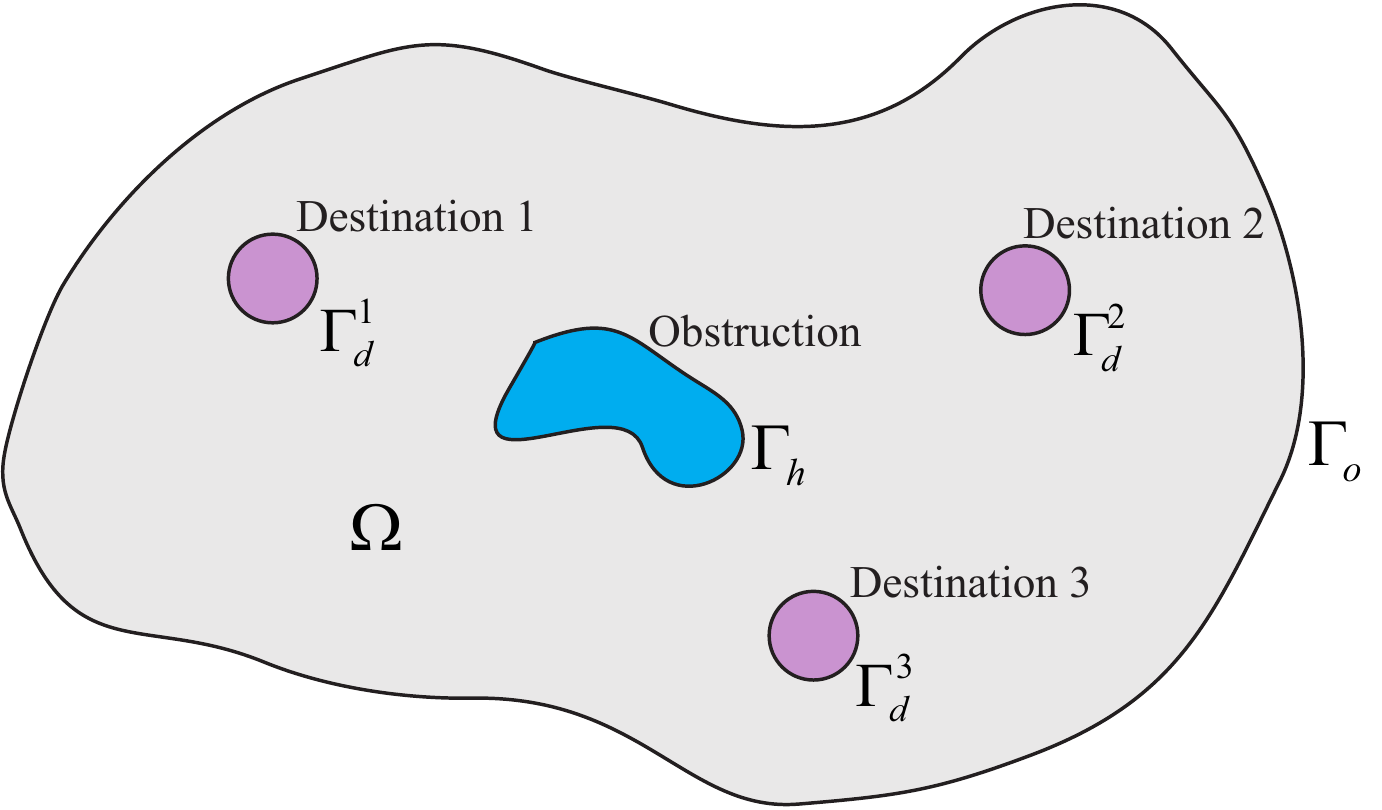}
    \caption{Illustration of continuum space with three destination areas and one obstruction (lake).}
    \label{f1}
\end{center}
\end{figure}

As mentioned earlier, continuum pedestrian models and continuum-space DTA models are  mathematically similar.  The first formulation of this type of models was proposed in the context of pedestrian dynamics by \textbf{\citet{Hughes2002Continuum}}.  
In this study, a two-dimensional walking infrastructure is represented as a continuum with domain $\Omega \subset \mathbb{R}^2$, as seen in Fig. \ref{f1}. $\Gamma_o$ denotes the outer spatial boundary, $\Gamma_h$ is the hard boundary of any obstruction at which no traveler is allowed to enter or exit the walking facility, and $\Gamma_d$ represents the boundary of the destination area. In more general cases, where there are more than one destination area in the study, $\Gamma_d^i$ will represent the boundary of each destination area, $i$.

Let 
$\rho(x,y,t)$ 
be the time-varying density of travelers, at location $(x, y)$ at time $t$,
$\bf{u}$$=(u_1(x, y, t),\ u_2(x, y, t))$ be the velocity vector, 
and
$\mathbf{f}(x,y,t) = (f_1(x, y, t),\ f_2(x, y, t)) \equiv \rho \mathbf{u}$ be the flux vector.
The conservation of travelers yields
\begin{equation} \label{eqHughesConservation}
\rho_t+\nabla \cdot \mathbf{f}=0, \hskip 4.6cm \mbox{(2D conservation law)}
\end{equation}
where $\nabla \cdot \mathbf{f}=\partial f_1 / \partial x+ \partial f_2 / \partial y$ is the (spatial) divergence of $\mathbf{f}$.
Hughes makes three hypotheses in order to model the nature of the pedestrian flow. The first hypothesis asserts that the speed of pedestrians at each point is an isotropic function of the density at that point and time
\begin{equation} \label{eqHughesSpeedByDensity}
\norm{\mathbf{u}}=V(\rho), \hskip 6.2cm \mbox{(isotropic case)}
\end{equation}
where $\ \norm{\mathbf{u}} \equiv \sqrt{u_1^2+u_2^2}$ is the norm of the velocity vector and $V(\rho)$ represents the speed-density relationships (and plays a role of the MFD in this review). 

In the second hypothesis, \citeauthor{Hughes2002Continuum} defines a potential function, $\phi(x,y,t)$, that represents the cost of reaching the destination starting from $(x,y,t)$ along the minimum cost path. Therefore, the  motion of any pedestrian is in the direction with maximum potential reduction, i.e. in the direction perpendicular to the isopotential curves. 
Since the gradient vector of the potential function, $\nabla \phi \equiv (\frac{\partial \phi}{\partial x},\frac{\partial \phi}{\partial y})$, is perpendicular to these isopotential curves in the direction of maximum increase, it follows that the flux vector, and as a result the velocity vector, are parallel to the gradient of the cost potential and in the \textit{opposite} direction, i.e.:
\begin{equation} \label{eqHuangRouteChoiceB}
\mathbf{u} \ / \! / \ \mathbf{f} \ / \! / \ -\nabla \phi. \hspace{4.9cm} \mbox{(optimum movement direction)}
\end{equation}
where // means ``is parallel to". Thus, the velocity vector is given by
\begin{equation} \label{eqHughesVelocityByDensity}
\mathbf{u}=-\frac{\nabla \phi}{\norm{\nabla \phi}} \ V(\rho), \hskip 3.7cm \mbox{(speed vector parallel to streamlines)}
\end{equation}
where $\nabla \phi / \norm{\nabla \phi}$ is the unit gradient vector of the potential function. 

The third hypothesis defines a discomfort function, $g(\rho)$, to account for the behavior of the pedestrians to avoid higher densities, which satisfies 
\begin{equation} \label{eqHughesDiscomfort}
gV=1/||\nabla \phi||, \hskip 3.9cm \mbox{(discomfort function $g$)}
\end{equation}
where $g \geqslant 1, \partial g/\partial \rho\geqslant 0 .$
However, \citet{HUANG2009Revisiting} explain that assuming the discomfort factor as an increasing function of density makes pedestrians perceive faster movement at higher densities, which is counterintuitive. Instead, they propose $\partial g/\partial \rho\leqslant 0$.

With all, Hughes' model can be expressed as
\begin{subequations} \label{Hughes}
    \begin{empheq}[left = \empheqlbrace\,]{align}
      & \rho_t-\nabla \cdot \left(\rho \, g\, V^2\, \nabla \phi\right)=0,  & \hspace{0.5cm} \mbox{(Hughes' model)} \label{HughesA}\\
       & gV=1/||\nabla \phi||, \label{eqTransform1B}
    \end{empheq}
\end{subequations}
which is not a hyperbolic system, as customary in the 1D traffic flow literature.

\citeauthor{Hughes2002Continuum} extends his model to incorporate multiple pedestrian types to represent different walking characteristics and destinations. Segregated potential functions are defined for each pedestrian type and the same hypotheses are developed for different pedestrian types.
The most important point in this extension is that the speed of a pedestrian type at each point and time is determined by the total density of all pedestrian types rather than the density of a single pedestrian type.

\citeauthor{Hughes2002Continuum}' paper 
has been the seminal work in the area of dynamic traffic assignment using continuum space, for both pedestrian and vehicular traffic. The following subsections will review the existing literature in this topic.

\subsection{Reactive Dynamic User Equilibrium Models}
In the Reactive Dynamic User Equilibrium (RDUE) problem, travelers  choose the route  that minimizes the instantaneous travel cost and change their choice en route as a result \citep{Boyce1995Solving}. \textbf{\citet{HUANG2009Revisiting}} revisited \citeauthor{Hughes2002Continuum}' \cite{Hughes2002Continuum} model and demonstrated that the route choice strategy in the model satisfies the RDUE principle. To see this, let $c(x, y, t)$ be the local cost per unit distance of movement,
\begin{equation} \label{eqHuangRouteChoiceC}
c=\norm{\nabla \phi},
\end{equation}
which is an Eikonal equation. They show that to minimize the instantaneous walking cost of a traveler, the route choice strategy should satisfy 
\begin{equation} \label{eqHuangRouteChoiceA}
c \ \frac{\mathbf{f}}{\norm{\mathbf{f}}}+\nabla \phi=0. \hspace{3cm} \mbox{(reactive route-choice strategy)}
\end{equation}
They prove that if $\phi$ is set to zero for $(x,y) \in \Gamma_d$, the potential function can be interpreted as the minimum instantaneous travel cost to the destination. 
They propose a general modeling structure consisting a set of partial differential equations as
\begin{subequations} \label{eqHuangPDE}
    \begin{empheq}[left = \empheqlbrace\,]{align}
       \rho_t+\nabla \cdot \mathbf{f}= 0, 
      \label{eqHuangPDEa}\\
       \norm{\mathbf{f}}= \rho\, V, 
      \label{eqHuangPDEb}\\
        c\frac{\mathbf{f}}{\norm{\mathbf{f}}}+\nabla \phi=0, 
       \label{eqHuangPDEc}
    \end{empheq}
\end{subequations}
subject to appropriate initial and boundary conditions, typically:
\begin{subequations} \label{eqHuangInitial}
    \begin{empheq}[left = \empheqlbrace\,]{align}
    	& \rho(x,y,0)=\rho_0(x,y), & \forall(x, y) \in \Omega, \qquad \label{eqHuangInitialA}\\
      	& \mathbf{f} \cdot \mathbf{\hat{n}}(x,y)=q(x,y,t), & \forall(x, y) \in \Gamma_o, \qquad \label{eqHuangInitialB}\\ 
       & \phi=0, & \forall(x, y) \in \Gamma_d, \qquad \label{eqHuangInitialC}
    \end{empheq}
\end{subequations}
where $\mathbf{\hat{n}}(x,y)$ is the unit normal vector toward $\Omega$ on $\Gamma_o$, $q(x,y,t)$ denotes the demand flow crossing $\Gamma_o$ and $\rho_0(x,y)$ is the initial density inside $\Omega$.

\textbf{\citet{Jiang2011Macroscopic}} extend their previous pedestrian flow \citep{Jiang2009Reactive}  by introducing a non-decreasing discomfort function and find the direction of movement of pedestrians, $\hat{d}$, using
\begin{equation} \label{eqJiang2011MovementDirection}
\hat{d}=\frac{\mathbf{u}}{\norm{\mathbf{u}}}=\frac{\mathbf{f}}{\norm{\mathbf{f}}}=\frac{\mathbf{w_1}+\epsilon \mathbf{w_2}}{\norm{\mathbf{w_1}+\epsilon \mathbf{w_2}}}, \qquad \qquad \qquad  \epsilon \geqslant 0,
\end{equation}
where vectors $\mathbf{w_1}=-\frac{\nabla \phi}{\norm{\nabla \phi}}$ and $\mathbf{w_2}=-\frac{\nabla g}{\sqrt{1+\norm{\nabla g}^2}}$ are directed along the dynamic paths of the shortest instantaneous travel time and the higher comfort level, respectively, and $\epsilon$ is a positive constant representing the psychological influence of avoiding higher densities. They reconstruct the pedestrian flow model \eqref{eqHuangPDE} by replacing   \eqref{eqHuangPDEc} with   \eqref{eqJiang2011MovementDirection} subject to the initial and boundary conditions of   \eqref{eqHuangInitial} and adding
\begin{equation} \label{eqJiang2011Initial}
\mathbf{f} \cdot \mathbf{\hat{n}}=0, \qquad \qquad \qquad \qquad \qquad \qquad \quad \forall(x,y) \in \Gamma_h,
\end{equation}
to reflect that no traveler enters or exits the obstructions inside the continuum domain.
\textbf{\citet{Jiang2014Macroscopic}} use RDUE in continuum space for dynamic traffic assignments of vehicles inside an urban network. 
The urban network is considered as a 2D continuum with a configuration similar to Fig. \ref{f1}, where the destination areas are the central business districts (CBD) inside the urban network. The local travel cost per unit distance is reformulated as
\begin{equation} \label{eqJiang2014LocalCost}
c(x,y,t)=\kappa \left( \frac{1}{U}+\pi_0\rho^2 \right), \hspace{1.5cm} \mbox{(general local cost function)}
\end{equation}
where $\kappa$ is the value of time and $\pi_0$ is a nonnegative constant interpreting other costs, such as avoiding higher densities. Their model is the following conservation law (CL):
\begin{subequations} \label{eqJiang2014PDE}
    \begin{empheq}[left = CL: \empheqlbrace\,]{align}
      & \rho_t+\nabla \cdot \mathbf{f} = q, & \forall(x, y) \in \Omega, \qquad \label{eqJiang2014PDEa}\\
      & \mathbf{f}=-\norm{\mathbf{f}}\frac{\nabla \phi}{\norm{\nabla \phi}} & \forall(x, y) \in \Omega, \qquad \label{eqJiang2014PDEb}\\
       & c=\norm{\nabla \phi}, & \forall(x, y) \in \Omega, \qquad \label{eqJiang2014PDEc}
    \end{empheq}
\end{subequations}
subject to initial and boundary conditions similar to   \eqref{eqHuangInitial}, where   \eqref{eqHuangInitialB} is replaced by   \eqref{eqJiang2011Initial}.

\textbf{\citet{Jiang2009Reactive}} develop a bi-directional pedestrian flow model, where two pedestrian groups, $a$ and $b$, have local walking speeds $U^a$ and $U^b$ dependent on the densities and direction of movement of both pedestrian groups, i.e.:
\vspace{-10pt}
\begin{subequations} \label{eqJiang2009Speed}
    \begin{empheq}[left = \empheqlbrace\,] {align}
       & U^a(\rho^a,\rho^b,\Psi)=U_f\ e^{-\alpha\left(\rho^a+\rho^b\right)^2}\ e^{-\beta (1-\cos(\Psi))(\rho^b)^2}, \label{eqJiang2009SpeedA}\\
       & U^b(\rho^a,\rho^b,\Psi)=U_f\ e^{-\alpha\left(\rho^a+\rho^b\right)^2}\ e^{-\beta (1-\cos(\Psi))(\rho^a)^2}, \label{eqJiang2009SpeedB}
    \end{empheq}
\end{subequations}
where $\rho^a(x,y,t)$ and $\rho^b(x,y,t)$ are densities of each group, $\Psi(x,y,t)$ is the intersecting angle between two pedestrian groups, $U_f$ is the free-flow walking speed of pedestrians, and $\alpha$ and $\beta$ are the model parameters. They also assume that the walking cost is equal to the travel time. Using these definitions and assumption they propose segregated systems of PDEs and initial and boundary conditions similar to   \eqref{eqHuangPDE} and \eqref{eqHuangInitial}, respectively, for each pedestrian group. 
They observed that pedestrians walking to the same destination form stripes, which shows the cooperative strategy of pedestrians moving in different directions, as discussed in \citet{Dammen2003Experimental}.

\textbf{\citet{Xiong2011High}} propose a high-order computational scheme for the numerical solution of this problem.
They found that the high-order scheme can obtain the same resolution on much coarser meshes with less computation time, compared with the lower-order schemes.

\textbf{\citet{Huang2009Dynamic}} extend the bi-directional pedestrian flow model developed by \citet{Jiang2009Reactive} by introducing look-ahead behavior in order to induce a viscosity effect on movement patterns of pedestrians. They first define the look-ahead displacement as
\begin{equation} \label{eqHuang2009Displacement}
(\delta x^i, \delta y^i)=\epsilon\ \frac{\mathbf{f^i}}{\norm{\mathbf{f^i}}}, \hspace{5cm} i \in \{a,b\},
\end{equation}
where $\epsilon>0$ is a small look-ahead distance of a pedestrian following the crowd in front of him corresponding to the traffic condition ahead at location $(x+\delta x^i,y+\delta y^i)$. The walking speed function with look-ahead effect for each group of pedestrians is given by
\begin{subequations} \label{eqHuang2009Speed}
    \begin{empheq}[left = \empheqlbrace\,] {align}
       & \ddot{U}^a(x,y,t)=U^a \left(\rho^a(x+\delta x^a,y+\delta y^a,t),\rho^b(x,y,t), \Psi \right), \label{eqHuang2009SpeedA}\\
       & \ddot{U}^b(x,y,t)=U^b \left(\rho^a(x,y,t),\rho^b(x+\delta x^b,y+\delta y^b,t), \Psi \right), \label{eqHuang2009SpeedB}
    \end{empheq}
\end{subequations}
where $U^a$ and $U^b$ are as defined in   \eqref{eqJiang2009Speed}. The rest of the model is the same as the model developed in \citet{Jiang2009Reactive}.
They proposed a numerical solution algorithm for the model, but they did not calibrate or validate this solution algorithm.


\textbf{\citet{Jiang2012Numerical}} revisit their previous work \citep{Jiang2009Reactive} for bi-directional pedestrian flow and introduce a more general local cost function, rather than assuming it to be equal to the travel time, as
\begin{equation} \label{eqJiang2012LocalCost}
c^i(x,y,t)=\frac{1}{U^i(\rho^a,\rho^b,\Psi)}+\gamma\, g^i(\rho^a,\rho^b,x,y,t), \hspace{1.3cm} i \in \{a,b\},
\end{equation}
where $\gamma$ is a weighted factor and $g^i$ represents the discomfort function for pedestrian group $i$. Then, they reconstruct their previous model by substituting the local cost function with   \eqref{eqJiang2012LocalCost}. 
They claim that numerical results of the model conform well to the experimental data and show that the model is rational and efficient.

\subsection{Predictive Route Choice}
In the Predictive Dynamic User Equilibrium (PDUE) problem travelers choose their routes to minimize the actual  travel cost they will experience, in contrast to the RDUE problem, where travelers seek to minimize their instantaneous travel cost. Thus, the cost potential function $\phi(x,y,t)$ in equilibrium is now the \textit{actual} travel cost (to the destination incurred by a traveler who leaves location $(x,y)$ at time $t$ toward the destination). 
\textbf{\citet{JIANG2011Dynamic}} propose a PDUE model for vehicular traffic for an urban network with only one CBD. They define an isotropic speed function as 
\begin{equation} \label{eqJiang2011bSpeed}
V(x,y,\rho)=U_f(x,y)e^{-\beta \rho^2},\hspace{0.7cm} \mbox{(isotropic speed-density relationship)}
\end{equation}
where $\beta$ is a positive parameter influenced by the road condition and etc. The local cost per unit distance is given by
\begin{equation} \label{eqJiang2011bLocalCost}
c(\rho)=\kappa \left( \frac{1}{V(\rho)}+\pi\left( \rho \right) \right), \hspace{1.5cm} \mbox{(general local cost function)}
\end{equation}
where $\pi(\rho)$ is a function representing other density-dependent costs.
They assert that the PDUE principle is satisfied, if the following conditions are met
\begin{subequations} \label{eqJiang2011bPDUEcondition}
    \begin{align} 
       & (u_1,u_2,1) \parallel - \bar{\nabla} \phi
       , \label{eqJiang2011bPDUEconditionA}\\
       &\norm{\bar{\nabla}\phi} =\frac{cV}{\sqrt{V^2+1}}
       , \label{eqJiang2011bPDUEconditionB}
    \end{align}
\end{subequations}
where $\overline{\nabla}\phi=\left( \frac{\partial \phi}{\partial x},\frac{\partial \phi}{\partial y}, \frac{\partial \phi}{\partial t}\right)$ is the potential gradient including the time component. Then, they develop a vehicular PDUE model based on these definitions and assumptions.
Nevertheless, \textbf{\citet{Du2013Revisiting}} revisit this model and claim that the path-choice strategy is not appropriate, e.g., when the cost function is time-independent, $\phi_t=0$, the condition stated by   \eqref{eqJiang2011bPDUEconditionA} cannot be satisfied.
They also indicate that there is an inconsistency in the units of one the PDE equations in the model and they prove that this equation is redundant. 

Their final model consists of a conservation law (CL) part similar to   \eqref{eqJiang2014PDE} and a Hamilton-Jacobi (HJ) part, governing $\phi(x,y,t)$:
\begin{subequations} \label{eqDu2013HJ}
    \begin{empheq}[left = HJ:\empheqlbrace\,]{align}
      & \frac{1}{V}\phi_t-\norm{\nabla \phi}=-c, & \forall(x, y) \in \Omega, \qquad \label{eqDu2013HJa}\\
       & \phi(x, y, t)=\phi_{CBD}, & \forall(x, y) \in \Gamma_d, \qquad \label{eqDu2013HJb}\\
       & \phi(x, y, t_{end})=\phi_0(x,y), & \forall(x, y) \in \Omega. \qquad \label{eqDu2013HJc}
    \end{empheq}
\end{subequations}

They interpreted $\phi_{CBD}$ as the cost of entering to the CBD and $\phi_0(x,y)$ as the instantaneous travel cost from any point at time $t=t_{end}$ to the CBD, when there is no traffic in the city and the travel cost is only related to the travel time, computed by an Eikonal equation as
\begin{subequations}
    \begin{empheq}[left = \empheqlbrace\,]{align}
      & \norm{\nabla \phi_0}=c(x, y, t_{end}), &\qquad \forall(x, y). \in \Omega \qquad \label{eqEikonala}\\
       & \phi_0=\phi_{CBD}, &\qquad \forall(x, y) \in \Gamma_d. \qquad \label{eqEikonalb}
    \end{empheq}
\end{subequations}

Notice that the initial time in the CL portion is $t=0$, whereas in the HJ portion the initial time is $t=t_{end}$.
\textbf{\citet{Lin2017Predictive}} extend the model proposed by \citet{Du2013Revisiting} to a polycentric urban city, where the speed and local cost per unit distance are computed separately for each  group of vehicles based on   \eqref{eqJiang2011bSpeed} and \eqref{eqJiang2011bLocalCost}, respectively, using the total density instead of the individual density of each group.
\textbf{\citet{Hoogendoorn2004Dynamic}} develop a pedestrian PDUE model for multiple pedestrian groups, where each group has a different destination, as
\begin{subequations} \label{eqHoogendoornPDE}
    \begin{empheq}[left = \empheqlbrace\,] {align}
       & \rho^i_t+\nabla \cdot \mathbf{f^i}=q^i(x,y,t),  \hspace{3.8cm} \mbox{(CL equation)}\label{eqHoogendoornPDEb}\\
       & -\phi^i_t=\min_{\mathbf{\tilde{u}}}\{\mathcal{L}(x,y,t,\mathbf{\tilde{u}})+\mathbf{\tilde{u}}\cdot\nabla\phi^i\}, \hspace{1.2cm} \mbox{(HJB equation)} \label{eqHoogendoornPDEa}\\
       & \mathbf{u^i}=arg \min\{\mathcal{L}(x,y,t,\mathbf{\tilde{u}}) +\mathbf{\tilde{u}}\cdot\nabla\phi^i\}, \label{eqHoogendoornPDEc}
    \end{empheq}
\end{subequations}
where, $\mathcal{L}$ is the running cost, $\mathbf{\tilde{u}}$ is any possible velocity vector and $\mathbf{u^i}$ is the optimal velocity vector of group $i$. This model is comprised of two coupled PDEs; the first PDE is a Hamilton-Jacobi-Bellman (HJB) equation and the second PDE is a CL equation. In order to solve the model numerically, the authors presented a heuristic iterative approach, where the HJB and the CL equations are calculated in each iteration, respectively. Thus, the minimum value problem has to be solved twice in each iteration. 

\textbf{\citet{Du2015Reformulating}} revisit \citeauthor{Hoogendoorn2004Dynamic}'s (\citeyear{Hoogendoorn2004Dynamic}) model for vehicular traffic for the sake of simplifying the numerical solution. They reduce the HJB equation to a HJ equation by removing the minimum value problem from the solution procedure and computing the actual travel cost when the velocity vector is known. Hence, the minimum value problem is only solved once in each iteration in the CL portion of the model. They first reconsider the speed function as an anisotropic function as
\begin{equation} \label{eqDu2015Speed}
U^i=U_f\ h \left( \psi^i(x,y,t) \right) \ g\left( \sum_{i} \rho^i \right),\hspace{0.2cm} \small \mbox{(anisotropic speed-density function)}
\end{equation}
where $U_f$ is the free-flow speed, $\psi^i$ is the angle between the direction of movement of group $i$ and the x-axis, $h(\psi^i)$ denotes the adjustment to the free-flow speed due to the direction of movement, and $g(\sum \rho^i)$ is a monotonically decreasing discomfort function. The local travel cost per unit distance for each group is defined as   \eqref{eqJiang2011bLocalCost}.
The proposed model consists of two systems of PDEs: the CL portion and the HJ portion. 
The HJ portion of the model for each group is the same as   \eqref{eqDu2013HJ}
and the CL PDEs and their initial and boundary conditions are expressed as:
\begin{subequations} \label{eqDu2015CLpde}
    \begin{empheq}[left =\empheqlbrace\,] {align}
       & \rho^i_t+\nabla \cdot \mathbf{f^i}=q^i, & \forall (x,y) \in \Omega, \label{eqDu2015CLpdeA}\\
       & U^i=U_f\, h\, g, & \forall (x,y) \in \Omega, \label{eqDu2015CLpdeB}\\
       & \psi^i =arg \min_{\tilde{\psi}}{p^i(x,y,t,\tilde{\psi})}, & \forall (x,y) \in \Omega, \label{eqDu2015CLpdeC}\\
       & \mathbf{f^i}\cdot \mathbf{\hat{n}}=0, & \forall (x,y) \in \Gamma_o \cup (\cup_{m\neq i}) \Gamma_d^m, \label{eqDu2015CLpdeD}\\
       & \rho^i(x,y,0)=\rho_0, & \forall (x,y) \in \Omega, \label{eqDu2015CLpdeE}
    \end{empheq}
\end{subequations}
where the CBDs other than the destination CBD of group $i$ are viewed as obstructions for the travelers of group $i$ and $p^i$ is an auxiliary function defined as
\begin{equation} \label{eqDu2015P}
p^i(x,y,t,\tilde{\psi})\equiv \nabla \phi^i(x,y,t) \cdot \mathbf{\tilde{u}^i} (x,y,t,\tilde{\psi}) +\tilde{c}^i(\tilde{\psi}) \tilde{U}^i(\tilde{\psi}).
\end{equation}

\subsection{Combination of Reactive and Predictive Route Choice}
We see that the reactive and predictive route choice behaviors have been combined in some papers for more realistic models. \textbf{\citet{Xia2009Dynamic}} propose a pedestrian DUE model founded on the hypothesis that "the pedestrians seek to minimize their estimated travel cost based on memory but temper this behavior to avoid high densities". In this hypothesis, the memory effect is the predictive or global portion and the behavior to avoid higher densities is the reactive or local portion of the route choice strategy. The route choice strategy is stated as
\begin{equation}\label{eqXiaRouteChoice}
\mathbf{f}(x,y,t) / \! / -\nabla \phi(x,y) - \omega \nabla c\left(\rho\right),
\end{equation}
where $\phi(x,y)$ is the time-independent minimum travel cost based on memory, $c(\rho)$ denotes the density-dependent costs per unit distance of movement at time $t$, and $\omega$ is a positive constant representing the psychological influence. The function $\phi(x,y)$ is given by an Eikonal equation and the function $c(\rho)$ is given by 
\begin{equation}\label{eqXiaDensityCost}
c(\rho)=\frac{1}{V(\rho)}+\beta\, g(\rho),
\end{equation}
where $V(\rho)=U_{f}\left(1-\frac{\rho}{\rho_{max}}\right)$ is the walking speed, $g(\rho)=\rho^2$ is the discomfort function, $\beta$ reflects the sensitivity of the pedestrians' route choice to discomfort, and $\rho_{max}$ is the jam density. 
\textbf{\citet{Hoogendoorn2015Continuum}} suggest another model for combining the global and local route choice behaviors, where they consider multiple groups of pedestrians with different destinations. The composite cost function for each group, $i$, is defined as
\begin{equation}\label{eqHoogendoorn2015Composite}
\omega_i=\phi_i+\varphi_i,
\end{equation}
where $\phi_i$ is the global cost function representing the minimum cost of getting to the destination, $\Gamma_d^i$, from $(x,y,t)$, and $\varphi_i$ is the local cost function reflecting additional walking costs due to local (unforeseen) fluctuations in the density, which is composed of crowdedness and delay components, as 
\begin{equation}\label{eqHoogendoorn2015Local}
\varphi_i=\varphi_i^{crowdedness}+\alpha_i \varphi_i^{delay},
\end{equation}
where $\alpha_i$ is the weight representing the relative importance of the two factors for each pedestrian group. The route choice strategy gives the direction of movement as
\begin{equation}\label{eqHoogendoorn2015Direction}
\mathbf{f_i} / \!  / -\nabla \omega_i = -\nabla \phi_i - \nabla \varphi_i.
\end{equation}

The crowdedness term in local cost function describes the tendency of the pedestrians to avoid areas with higher densities and the delay term reflects the expected increase in delay caused by local densities. The authors also did some simulation experiments in order to investigate the impact of different factors introduced in their new model.

\subsection{Higher-order macroscopic models}
The first-order models that have been discussed here, are based on the assumption that the traffic flow is always in the equilibrium state and are not capable of describing non-equilibrium phenomena such as stop-and-go waves and formation of congestions. 
\textbf{\citet{Jiang2010Higher}} introduce a second-order pedestrian flow model incorporating an equilibrium of linear momentum equation in addition to the conservation of mass   \eqref{eqHughesConservation}, as
\begin{equation}\label{eqJiang2010Momentum}
\mathbf{u_t}+(\mathbf{u} \cdot \nabla) \mathbf{u}=\frac{V(\rho)\, \hat{d}-\mathbf{u}}{\tau} -R'(\rho)\, \frac{\nabla \rho}{\rho}, \hspace{0.5cm} \small \mbox{(second-order conservation law)}
\end{equation}
where $\mathbf{u}(x,y,t)$ is the average pedestrian velocity vector, 
$\hat{d}=(d_1,d_2)$ denotes the unit optimal movement direction, $\tau$ is the relaxation time of $\mathbf{u}$ toward the optimal velocity (taken as 0.5 s by the authors), $R(\rho)$ denotes the traffic pressure, which describes the response of pedestrians to compression, and $R'(\rho)=c_s^2(\rho)$, where $c_s$ is the sonic speed at which small disturbances propagate in a crowd. The first term at the right side of   \eqref{eqJiang2010Momentum} is called the relaxation term, which represents the relaxation to equilibrium, and the second term is called the anticipation term, which reflects the reaction to the surrounding traffic condition.   \eqref{eqJiang2010Momentum} can be rewritten more explicitly in the components form as
\begin{subequations} \label{eqJiang2010Components}
    \begin{empheq}[left = \empheqlbrace\,]{align}
      & \frac{\partial u_1}{\partial t}+u_1 \frac{\partial u_1}{\partial x}+u_2 \frac{\partial u_1}{\partial y}= \frac{V\, d_1-u_1}{\tau} -\frac{c_s^2}{\rho} \frac{\partial \rho}{\partial x}, \label{eqJiang2010ComponentX}\\       
       & \frac{\partial u_2}{\partial t}+u_1 \frac{\partial u_2}{\partial x}+u_2 \frac{\partial u_2}{\partial y}= \frac{V\, d_2-u_2}{\tau} -\frac{c_s^2}{\rho} \frac{\partial \rho}{\partial y}. \label{eqJiang2010ComponentsY}
    \end{empheq}
\end{subequations}

As usual,  the velocity vector  is assumed tangential to the gradient of $-\phi$, i.e. 
satisfying the Eikonal PDE  \eqref{eqHuangRouteChoiceC}.
The final RDUE model proposed in \citet{Jiang2010Higher} is a system of PDEs comprised of the CL   \eqref{eqHughesConservation}, the equilibrium of linear momentum   \eqref{eqJiang2010Momentum}, and the Eikonal  \eqref{eqHuangRouteChoiceC}. 
Using numerical experiments, they found that the traffic becomes more unstable as the anticipation factor decreases, which is explained more explicitly later in \citet{Jiang2016Comparison}.

\textbf{\citet{Jiang2016Macroscopic}} develop a second-order PDUE pedestrian flow model using 
\eqref{eqHughesConservation} and \eqref{eqJiang2010Momentum} for the CL part and 
\eqref{eqDu2013HJ} for the HJ part. 
They validated their model with  experimental pedestrian flow data collected under non-congested conditions. For congested pedestrian flow conditions, using numerical experiments, they found that increasing the anticipation factor results in significant reduction of the density near bottlenecks.

\textbf{\citet{Jiang2016Comparison}} compare the RDUE and PDUE models developed by \citet{Jiang2010Higher} and \citet{Jiang2016Macroscopic}, respectively, and evaluate their behavior by numerical experiments. In order to conduct these numerical experiments, the authors define the equivalent traffic sonic speed as
\begin{equation}\label{eqJiang2016Sonic}
c_s(\rho)=\sigma \rho,
\end{equation}
where $\sigma$ is the anticipation degree representing the anticipation behavior of pedestrians to compression. The RDUE and PDUE models are tested using different values for the anticipation degree and they found that if the pedestrians have strong anticipation consciousness, i.e. $\sigma > 0.33$, traffic flow is always in a stable state, less congestion occurs, and there are no significant difference between two models. On the contrary, if the anticipation degree falls below a critical value, e.g. $\sigma <0.33$, traffic instability occurs, density increases at certain points and results into congestion, and the RDUE model exhibits more uniform density distribution compared to the PDUE model.

\subsection{Dynamic System Optimum}

In dynamic system optimum (DSO) assignment problems the goal is minimizing the total travel cost for all of the travelers in the system. Solving DSO problem using macroscopic fundamental diagram is yet an underdeveloped area, with \textbf{\citet{Tao2014Dynamic}} being the only DSO model. 
This model is based on the assumption that the intelligent transportation system (ITS) and the advanced traveler information system (ATIS) have perfect information about the time-varying traffic status and the travelers choose and change their route, if needed, entirely on the information given by the ITS and ATIS. The local cost per unit distance is defined as
\begin{equation}\label{eqTaoLocal}
c(\rho)=\frac{\kappa}{V(\rho)},
\end{equation}
where $\kappa$ is the value of time and $V(\rho(x,y,t))\equiv U_f(x,y)e^{-\beta \rho}$. The 
system cost is expressed as
\begin{equation}\label{eqTaoSystemCost}
\Phi(\rho, \theta)=
\kappa\iiint_{\Omega \times T}  \rho \ dx\,dy\,dt, \hspace{1.4cm} \mbox{(system cost)}
\end{equation}
where $\theta (x,y,t)$ is the directed angle of speed, representing travelers' route choice. The DSO model is formulated as an optimization problem with a feasible region in the function space as
\begin{equation}\label{eqTaoMinimum}
min\ \ \Phi(\rho, \theta),
\end{equation}
subject to the constraints
\begin{subequations} \label{eqTaoConstraints}
\begin{empheq}[left = \empheqlbrace\,]{align}
&\rho_t+\nabla \cdot \mathbf{f}=q, & \forall (x,y) \in \Omega, \qquad \label{eqTaoConstraintA} \\
&V=\norm{\mathbf{u}}, & \forall (x,y) \in \Omega, \qquad \label{eqTaoConstraintB}\\
&\rho(x,y,0)=\rho(x,y,t_{end})=0, & \forall (x,y) \in \Omega, \qquad \label{eqTaoConstraintC}\\
&\rho(x,y,t)=0, & \forall (x,y) \in \Gamma_o \cup \Gamma_h. \qquad \label{eqTaoConstraintD}
\end{empheq}
\end{subequations}

Since the analytic properties of the proposed model were not studied by the authors, they were not able to develop a global optimal solution for the model. Instead, they derive a locally optimal solution with low computation cost. 
\subsection{One-dimensional Spatial model}

Unlike previous models, which are in two spatial dimensions, \textbf{\citet{Laval2017Minimal}} incorporate a single spatial dimension to describe the changes in the area covered by the MFD. They investigate DUE conditions on a single origin-destination pair with two alternative routes, a freeway with a fixed capacity $\mu_0$ and the surrounding cite-streets (CS) network, described by an o-MFD with capacity $\mu_1$. 
The authors first describe the single MFD dynamics and find analytical solution for special cases. Subsequently, they add a freeway alternative under two scenarios for the CS network: (i) constant network length, and (ii) variable network length to account for the spatial extent of congestion. In the first scenario, the combined system dynamics are given by the following ODE
\begin{subequations} \label{eqLavalUEode}
\begin{empheq}[left = \empheqlbrace\,]{align}
&k'(t)=\frac{b \varpi(t)-kv}{1-bmv'/v^2}, & \text{(reservoir dynamics)}, \qquad \label{eqLavalUEodeA} \\
&k(0)=k_0, & \text{(initial conditions)}, \qquad \label{eqLavalUEodeC}
\end{empheq}
\end{subequations}
where $k$ is the occupancy, $t$ is measured in units of free-flow travel time inside the MFD, $\tau ^*$,  $b$ is a positive constant defining the shape of the MFD, and $v(k)$ is a dimensionless version of  the  speed-occupancy MFD.  The dimensionless parameters $m$ and $\varpi$ are given by
\begin{subequations} \label{eqLavalParameters}
\begin{align}
&m=\frac{\mu_0}{\mu_1}, & \text{(Freeway to CS capacity ratio)}, \label{eqLavalM} \\
& \varpi(t) \equiv \frac{\lambda(t)-\mu_0}{\mu_1}=\frac{\lambda(t)}{\mu_1}-m, & \text{(MFD demand intensity)}. \label{eqLavalRho}
\end{align}
\end{subequations}

\begin{figure} [t]
\begin{center}
\includegraphics[width=120mm]{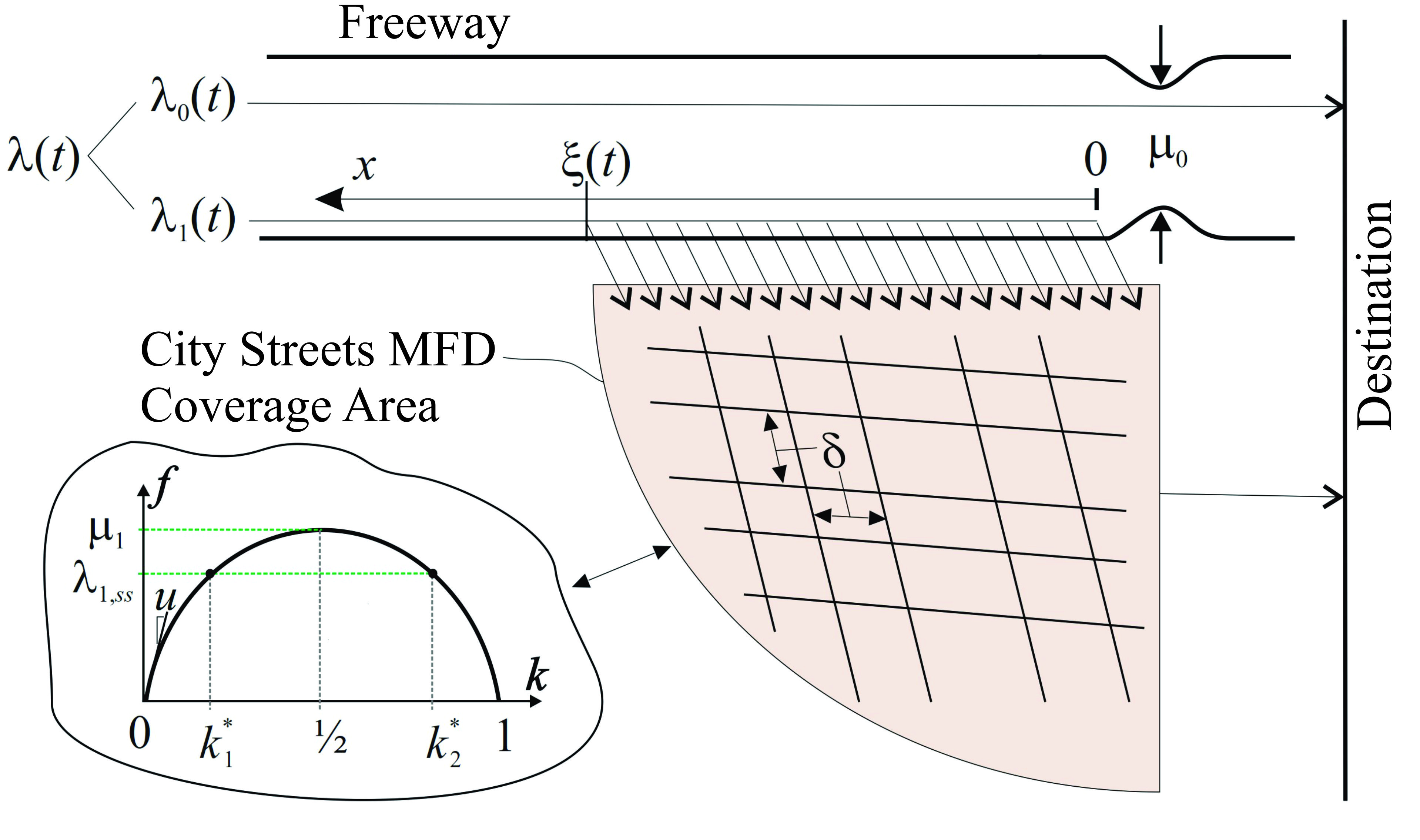}
\vspace{-10pt}
    \caption{Freeway vs. City streets network with continuum approximation for off-ramps \citep{Laval2017Minimal}.}
    \vspace{-10pt}
    \label{f2}
\end{center}
\end{figure}

The authors were  able to solve the ODE \eqref{eqLavalUEode} analytically for the autonomous case, where the demand is constant over time and the evolution of the system depends only on the occupancy. 
For scenario (ii) the authors incorporate continuum 
approximation (CA), proposed by \citet{Laval2009Graphical}, for off-ramps into their proposed formulation, as seen in Fig. \ref{f2}. The discrete off-ramps are considered as a continuum, where vehicles can exit the freeway to CS at any location upstream the bottleneck along $0 \leq x \leq \xi(t)$, where $x=0$ is the bottleneck location and $\xi(t)$ is a characteristic of DUE solution, named information wave, marking the most upstream location where vehicles divert from the freeway to the CS.

The implementation of CA in the proposed model makes it able to cope with time-varying network and trip lengths for the MFD. Notwithstanding the fact that analytical solutions become impossible for variable network and trip lengths, the numerical solutions strongly suggest that the proposed system have many analogies to the constant length model. The main difference is that gridlock does not happen in the proposed model and the steady-state solution is independent of surface network parameters, when time is expressed in units of MFD free-flow travel time.

\section{Solution Methods for Continuum-Space models}
\subsection{Analytical solutions for Hughes' model}
\begin{figure} [t]
\begin{center}
\includegraphics[width=110mm]{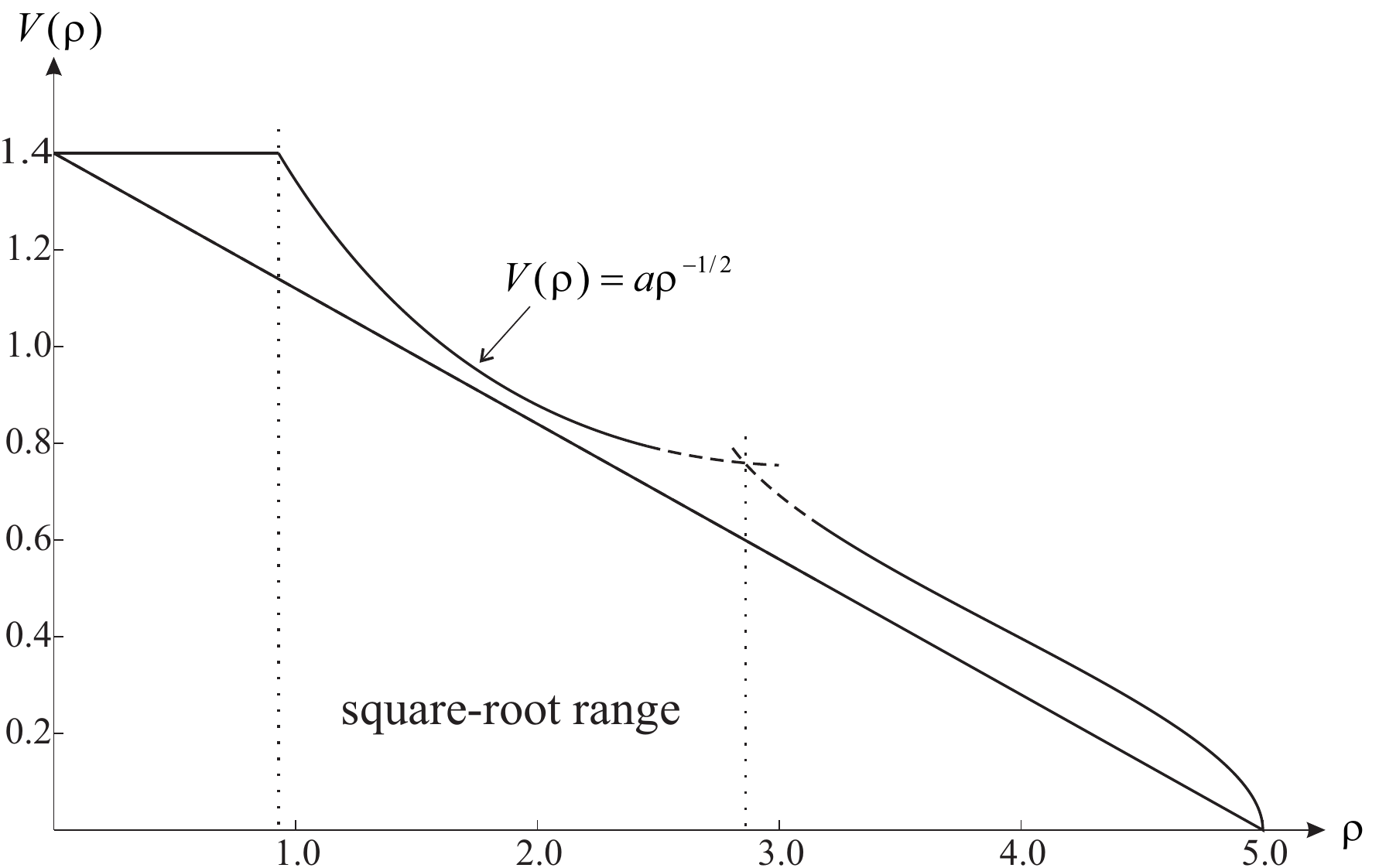}
    \caption{Square-root range in speed-density relationship \citep{Hughes2002Continuum}.}
    \label{f3}
\end{center}
\end{figure}

Hughes \cite{Hughes2002Continuum} shows that when the speed-density relation is proportional to $\rho^{-1/2}$, i.e.
\begin{equation}\label{eqSquareroot}
V(\rho)=a \rho^{-1/2},
\end{equation}
where $a $ is a positive constant, then his formulation has a symmetry under conformal mappings, i.e.   \eqref{eqHughesConservation} looks identical before and after the transformation.
These maps make it possible to obtain analytical solutions to problems with obstructions of complicated shape. There are map dictionaries where one can identify which particular mapping will transform the shape of the obstructions to lines or circles, for which the solution is known.

Since the square root approximation   (\ref{eqSquareroot}) is not a good representation of speed for all densities, Hughes assumes that it holds inside a density range, which we can call the ``square-root" range; see Fig. \ref{f3}. Unfortunately, when the density falls outside of this range, the symmetry is lost and the solution method becomes convoluted.

Later in the paper, we will show that the conformal mapping symmetry is valid for all density ranges and for any speed-density relationship. In fact,    \eqref{eqSquareroot} is needed only to simplify Hughes' model to Laplace's equation, 
\begin{equation} \label{eqLaplace}
  \nabla^2 \phi =0, \hspace{5cm} \mbox{(Laplace's equation)}
\end{equation}
whose solution is well understood. To see this notice that the conservation law   \eqref{HughesA} can be expressed as
\begin{equation} \label{eqTransform2}
\rho _t-\left(2 \rho  V'\nabla \rho  \cdot \nabla \phi  +(\rho  \nabla^2 \phi +\nabla \rho  \cdot \nabla \phi )V\right)V=0,
\end{equation}
by expanding the divergence term and assuming the discomfort function $g(\rho)=1$ without loss of generality. The reader can verify that the two terms of involving $\nabla \rho  \cdot
\nabla \phi $ cancel out under   \eqref{eqSquareroot} and since $\rho _t =0$ in steady-state, we end up with Laplace equation as sought.

Hughes also identifies two limit cases and proposes simpler equations that the solution satisfies in steady state: 

\vspace{-5pt}
\begin{enumerate}
\addtolength{\itemindent}{5mm}
\item at low  free-flow densities  $\nabla \rho  >>\nabla \phi $, therefore the CL   \eqref{HughesA} becomes
\begin{equation}
\nabla \rho  \cdot \nabla \phi \approx 0
\end{equation}
which means that lines of constant density are almost perpendicular to lines of constant potential.
\item at high congested densities  $\nabla \rho  << \nabla \phi $, thus the CL  becomes
\begin{equation}
\nabla^2 \phi \approx 0
\end{equation}
\end{enumerate}

Although not shown in the original paper, we can obtain these results using   \eqref{eqTransform2} and letting $\rho_t=0$ and $\nabla^2 \phi=0$ for the first case, and $\nabla \rho \approx 0$ for the second case.

\subsection{Numerical Solutions}
In order to solve continuum-space models, various numerical solutions are proposed in each study. Most of the solution procedures use a combination of standard solution methods for PDEs in order to discretize the continuous  space and time.
Conservation law and Hamilton-Jacobi equations are often solved by methods such as finite volume method (FVM), finite difference method (FDM), including Lax-Friedrichs (LF) and weighted essentially non-oscillatory (WENO) schemes, and finite element methods (FEM), including discontinuous Galerkin method (DGM).

Fast sweeping method (FSM) and fast marching method (FMM) are mostly utilized to solve the Eikonal equation. Most of the studies use total variation diminishing Runge-Kutta (TVDRK) for time discretization and/or integration.
Some of the PDUE studies have also viewed the model as a fixed-point problem and have solved the CL and HJ parts simultaneously using self-adaptive method of successive averages (SA-MSA).
Table \ref{t1} gives the different solution methods used in each study.

On a category on its own, 
\textbf{\citet{Hanseler2014Macroscopic}} develop a discrete-time discrete-space pedestrian flow model, named PedCTM, which utilizes a discretization scheme similar to Daganzo's Cell Transmission Model (CTM) \citep{Daganzo1994Cell,Daganzo1995Cell}. We have included this model in the continuum category since PedCTM can be viewed as a numerical solution method of a conservation law in 2D, similar to the CTM for one spatial dimension.

The period of analysis is discretized into intervals, $\tau$, with uniform length $\Delta t$. The 2D walking domain is partitioned into cells, $\varkappa$, with uniform shape, which are assumed to be partitioned orthogonally with spacing size of $\Delta L$. The pedestrians are assumed to be distributed homogeneously within a cell and their movements are not modeled explicitly.

A pedestrian group, $i$, is defined by a route $\eta_i$, a departure time interval $\tau_i$, and a size $X_i$. A route can be defined as a sequence of areas, $\eta=(r_o,r_1,...,r_d)$, where each area $r$ is a set of cells and $r_o$ and $r_d$ denote the origin and destination areas, respectively.
Inside each route, several sequences of cells can connect the origin area to the destination area, where each sequence of cells is called a path and the en-route path choice is made by computing turning proportions at each cell.

In order to describe the internal dynamics of the model, a normalized PedCTM model is developed, which is independent of the absolute values of free-flow speed and jam density.
They utilize an isotropic normalized cell-based version of the speed-density relationship proposed by \citet{Weidmann1992} as
\begin{equation} \label{eqHanselerFundamental}
\upsilon(\bar{\rho}_{\varkappa,\tau})=1-\exp \left(-\omega(\frac{1}{\bar{\rho}_{\varkappa,\tau}}-1)\right),\hspace{0.2cm} \small \mbox{(normalized speed-density function)}
\end{equation}
where $\upsilon$ is the normalized prevailing walking speed function, $\bar{\rho}_{\varkappa,\tau}$ represents the normalized density of cell $\varkappa$ at time interval $\tau$ and $\omega$ is a dimensionless parameter.
The route choice behavior of  pedestrians is described using route-specific cellular potential fields $P_{\varkappa,\tau}^\eta$, given by
\begin{equation} \label{eqHanselerPotential}
P_{\varkappa,\tau}^\eta=\alpha \digamma_\varkappa^\eta-\beta D_{\varkappa,\tau}, \hspace{4.6cm} \mbox{(potential field)}
\end{equation}
where $\digamma_\varkappa^\eta$ is the static floor field, which is assumed to be the minimum distance (number of cells) to be traversed from cell $\varkappa$ toward the destination area along route $\eta$, $D_{\varkappa,\tau}$ represents the dynamic floor field, which is taken as $\upsilon(k_{\varkappa,\tau})$ and $\alpha$ and $\beta$ are positive weights. If the set of all adjacent cells to $\varkappa$ that are part of route $\eta$ is denoted by $\Theta_\varkappa^\eta$, the turning proportion corresponding to link $\varkappa \to \chi \in \Theta_\varkappa^\eta$ is expressed as
\begin{equation} \label{eqHanselerTurning}
\delta_{\varkappa \to \chi, \tau}^\eta=\frac{\exp \left(-P_{\chi,\tau}^\eta\right)}{\sum_{\chi' \in \Theta_{\varkappa}^\eta} \exp \left(-P_{\chi',\tau}^\eta\right)}. \hspace{2cm} \mbox{(turning proportion)}
\end{equation}

In this manner, the route-specific potential field can be interpreted as disutility of each cell in the route toward the destination area and the pedestrians try to decrease their disutility along the route in a reactive manner. 
\textbf{\citet{Hanseler2017Dynamic}} propose an anisotropic discrete-time discrete-space pedestrian flow model based on a stream-based pedestrian fundamental diagram (SbFD). Unlike \citet{Hanseler2014Macroscopic} where areas were comprised of cells, in \citet{Hanseler2017Dynamic} each area, $\zeta$, contains a number of streams, $s$, and the set of streams associated with area $\zeta$ is denoted by $\Lambda_\zeta$.
Each route, $\eta$, consists of a pair of origin and destination nodes and a set of streams ,$\Lambda_\eta$, connecting them. The walking speed of stream $s$ in area $\zeta$ is computed by
\begin{equation} \label{eqHanseler2017Speed}
U_s=U_f \exp \left(-\vartheta (\frac{N_\zeta}{A_\zeta})^2\right) \prod_{s' \in \Lambda_\zeta} \exp \left({-\beta (1-\cos \theta_{s,s'}) \frac{M_{s'}}{A_\zeta}}\right),
\end{equation}
where $N_\zeta$ and $A_\zeta$ are the accumulation and surface size of area $\zeta$, $M_{s'}$ denotes the number of pedestrians in stream $\lambda'$, $\theta_{s,s'}$ is the intersection angle between streams $s$ and $s'$, and $\vartheta$ and $\beta$ are the model parameters. At each node $\Upsilon$, where the streams intersect, a potential value $P_{\Upsilon_s^d,\tau}^\eta$ is defined as the remaining walking time to destination $d$ using stream $s$ in route $\eta$ that can be calculated using any shortest path algorithm, e.g., \citet{Dijkstra1959Note}. A path choice strategy similar to \citet{Hanseler2014Macroscopic} is proposed, which incorporates a weighted logit-type model at each node.

The performance of the model is compared to a few isotropic specifications \citep{Wong2010Bidirectional,Weidmann1992} at the example of two case studies. The analysis reveals that the consideration of anisotropy improves the accuracy of the proposed model compared to the tested isotropic specifications. Overall, the model proposed in \citet{Hanseler2017Dynamic} can been seen as an analogy to the PedCTM model \citep{Hanseler2014Macroscopic}, where the cell-based fundamental diagram, potential fields and path choice are replaced by stream-based ones and anisotropy is also taken into account.

\section{Discrete-space models}

In discrete-space models, the modeling region is divided into a finite number of zones describing a tessellation of the $(x,y)$-plane. Each zone has a well-defined MFD with traffic dynamics given by the conservation ODE,   \eqref{model0}. 

A standard assumption in the literature
is that demand should be restricted by  the congested branch of the o-MFD, also called the supply function, $S(n) $:
\begin{equation}\label{capacity constraint}
	\lambda (t)\leqslant S (n (t)).  \hskip 4.5cm \mbox{(supply constraint)}
\end{equation}
to reflect that, in congestion, the inflow cannot exceed the outflow. This is consistent with the idea that all links in the region, including those on its perimeter, have the same (congested) flow and therefore if demand exceeds it, a queue accumulates outside the region.

\begin{figure} [t]
\begin{center}
\includegraphics[height=52mm]{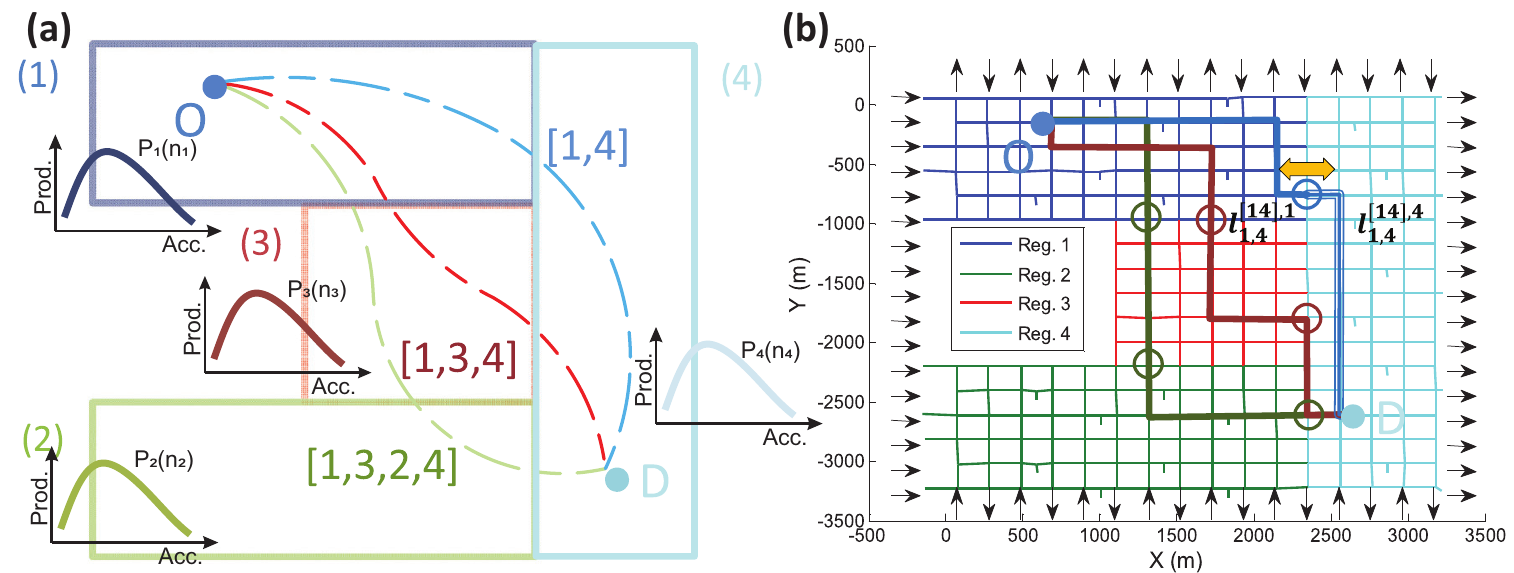}
\vspace{-5pt}
    \caption{Various paths from region 1 to region 4 shown in: (a) regional and (b) link-level representations  \citep{Yildirimoglu2014Approximating}.}
    \label{f4}
\end{center}
\end{figure}

\textbf{\citet{Yildirimoglu2014Approximating}} incorporate MFD into a DTA model in a heterogeneous urban network divided into homogeneous regions with low-scatter MFDs and find the  dynamic stochastic user equilibrium (DSUE) numerically. 
They define a path as the sequence of regions from the origin to the destination. Fig. \ref{f4}(a) shows three possible paths from region 1 to region 4 as sequences of different regions.
Since the trip lengths are time-dependent and there may be various link sequences complying with the path definition, stochastic network loading (SNL) is implemented to address variable trip lengths within and between regions. 

In order to determine the shortest path, the algorithm developed by \citet{Chabini1998Discrete} is implemented to determine the sequence of links leading to the minimum experienced travel time inside a sequence of regions using the link-level representation of the network, as seen in Fig. \ref{f4}(b),
while the speed on all links in a region is given by the average speed of the MFD, based on the accumulation level. 
Then, the probability of using each path is calculated using a multinomial logit discrete choice model based on the calculated minimum travel times.
Traffic equilibrium is formulated as a fixed-point problem that is solved in an iterative manner using method of successive averages (MSA), where the step size is determined a priori.

The proposed model is compared to a one-shot (non-iterative) assignment model with the aggregate dynamics, where the path choice is updated periodically in a way that the vehicles departing at each period choose their path based on instantaneous travel times measured at the previous period. 
It is demonstrated that the one-shot assignment model does not satisfy equilibrium conditions and overstate the network congestion, whereas the proposed model results in DSUE condition and less congestion.

\textbf{\citet{Yildirimoglu2015Equilibrium}} extend the work in \citet{Yildirimoglu2014Approximating} to a route guidance system based on dynamic system optimum conditions, where each region is divided into a number of sub-regions with well-defined MFDs and routes are redefined as 
sequences of sub-regions instead of sequences of links in \citet{Yildirimoglu2014Approximating}. Two macroscopic traffic models are developed in region-level and subregion-level, where the latter deploys a more detailed approach. The transfer of variables from the subregion-based model to the region-based model has also been described.

In this paper, DTA has been achieved using three different approaches: (i) by establishing PDUE conditions in subregion-level, (ii) by establishing DSO conditions in subregion-level, and (iii) by providing travelers with route guidance (RG) information based on DSO conditions in region-level, where the route guidance commands are applied in subregion-level. DSO conditions require equal and minimal marginal travel times on alternative routes at the same departure time. Several numerical studies have been conducted in order to assess different DTA approaches. The results reveal that RG approach outperforms DSO and PDUE approaches, except in very high demands which might be due to the inconsistency between subregion-level and region-level accumulation values near to the jam accumulation.

\begin{figure} [!b]
\begin{center}
\includegraphics[width=130mm]{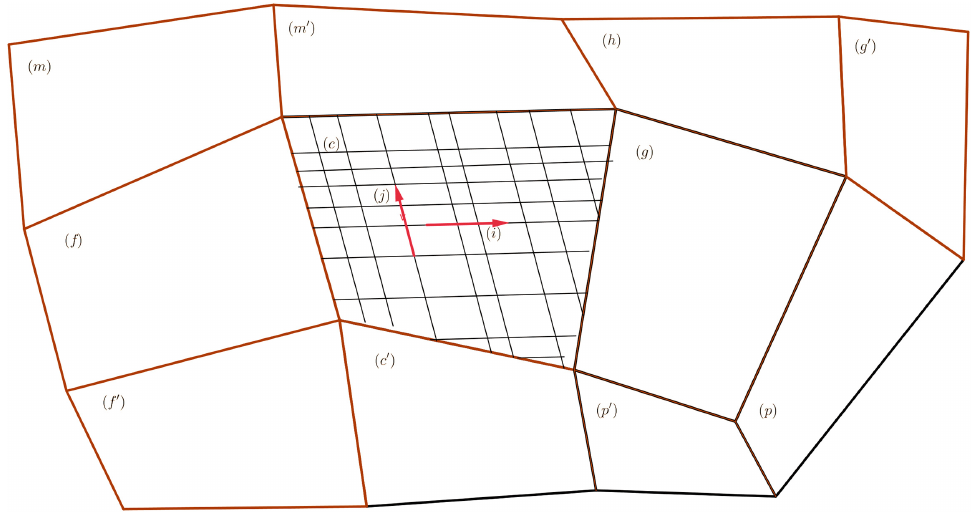}
    \caption{Network disaggregated into zones with quadrangular mesh of cells \citep{Sossoe2017Reactive}.}
    \label{f6}
\end{center}
\end{figure}

\textbf{\citet{Sossoe2017Reactive}} develop an RDUE model for bi-directional traffic flow based on the traffic flow model proposed in \citet{Saumtally2013Dynamical} and \citet{Sossoe2015Traffic}. The urban network is decomposed into zones, where each zone is meshed into quadrangular cells, as shown in Fig. \ref{f6}, and four inflow and outflow directions are considered for each cell. Cell flows and exit speeds in each direction are calculated by cell densities using intersection traffic flow model rules.

Based on the position of the target cell regarding the origin cell, only two possible outflow directions are assumed from the origin cell toward the designated target cell. Two possible paths are considered based on each of the possible directions and the turning proportions in each cell are computed using a weighted logit model, where the cost of each path is assumed to be the instantaneous travel time.

\subsection{Numerical Solutions}
In the case of discrete-space models, numerical methods are needed to solve (i) the single-MFD loading problem, and (ii) the network loading problem and equilibrium.

For the solution of (i), the main numerical component is the conservation ODE   \eqref{model0}, which is typically solved with  Euler's method, which is a first-order  finite differences scheme:
\begin{equation}\label{}
n(t+\Delta t)=n(t)+\Delta t(\lambda (t)-o (n))
\end{equation}
which converges linearly to the true solution as $\Delta t\rightarrow 0$.

\citet{Laval2017Minimal} derived analytical solutions for this problem for a family of demand curves well suited to model rush-hour periods, under the Greenshield (parabolic) MFD. These include  second degree polynomial, exponential, and logistic functions. For a general MFD, they found analytical solutions for the autonomous case, i.e. when the demand is time-independent.

In order to solve (ii), we have seen that the methods used are similar to the ones used in traditional DTA models with link exit-flow functions \citep{Merchant1978Model,Friesz1989Dynamic,Carey2004Exit}. Although the literature review on this vast subject is out of the scope of this manuscript, it is important to note that the lessons learned from those models can be directly applied to these simplified macroscopic networks.

\section{Conclusions and Outlook}

In this section, we summarize our main findings in terms of the observed trends in the literature, and identify issues  that deserve further research.
Table \ref{t1} presents an overview of the studied papers and their specifications. The lines between some rows delimit the paper reviewed in different sections and subsections of this paper.

\subsection{Differences between continuum-space pedestrian flow and vehicular DTA models}
We mentioned at the outset that  pedestrian and vehicular continuum DTA models are mathematically very similar, and therefore can be analyzed within the same framework. 
From a modeling perspective, however, there are some differences.
A quick look at Table \ref{t1} reveals that most of the studies investigating pedestrian flow incorporate reactive route-choice component, whereas the studies regarding vehicular flow subsume predictive route-choice component. 
This seems rational since it means that pedestrians usually do not have complete information about the conditions toward their destination and seek to reduce their instantaneous travel cost by choosing the shortest path and also avoiding higher densities, which cause discomfort and excessive delay. 
On the other hand, drivers are increasingly having better  information about the traffic conditions in different routes, using their own experience or traffic applications, and try to minimize their actual travel cost.

Another obvious difference between continuum-space pedestrian flow and vehicular  DTA models is that in the pedestrian flow models the demand is exogenous, i.e. set to zero inside the continuum domain and set to $q$ on outer boundary, while it is endogenous in vehicular DTA models, i.e. set to $q$ inside the continuum domain and set to zero on outer boundary.

\setlength{\abovecaptionskip}{-2pt}
\setlength{\belowcaptionskip}{10pt}

\begin{table}
\vspace{-60pt}
\setlength\tabcolsep{1pt}
\setlength\extrarowheight{0.5pt}
\begin{center}
\fontsize{7.5}{9}\selectfont
\centering
\rowcolors{1}{gray!25}{white}
\captionof{table}{Overview of studied papers}
\label{t1}
\begin{tabular} {K{2.1cm} K{1.45cm} K{1.7cm} K{1.4cm} K{1.4cm} K{1.65cm} K{1.65cm} K{5.3cm}}
\rowcolor{gray!75}
\textbf{Source} & \textbf{Granularity} & \textbf{Pedestrian or Vehicular} & \textbf{Route Choice} & \textbf{Multi-Destination} & \textbf{Directionality} & \textbf{Solution Methods} & \textbf{Main Contribution}\\
\hline
\citet{Hughes2002Continuum} & Continuum & Pedestrian & --- & Yes & Isotropic & Conformal maps & Develops the framework for continuum modeling of DTA \\ 
\citet{HUANG2009Revisiting} & Continuum & Pedestrian & Reactive & No & Isotropic & WENO, FSM, TVDRK &Proves that \citet{Hughes2002Continuum} satisfies the RDUE conditions \\ 
\citet{Jiang2011Macroscopic} & Continuum & Pedestrian & Reactive & No & Isotropic & DGM, TVDRK & Introduces a new route-choice strategy using  \eqref{eqJiang2011MovementDirection} \\ 
\citet{Jiang2014Macroscopic} & Continuum & Vehicular & Reactive & No & Isotropic & FVM, FSM, TVDRK & Extends the concept of RDUE to the vehicular traffic \\ 
\citet{Jiang2009Reactive} & Continuum & Pedestrian & Reactive & Yes, 2 & Bi-directional & FVM, FMM, TVDRK &  Extends the RDUE concept to bi-directional pedestrian flow \\ 
\citet{Xiong2011High} & Continuum & Pedestrian & Reactive & Yes, 2 & Bi-directional & WENO, FSM, TVDRK &  Proposes a different numerical solution procedure for the model developed in \citet{Jiang2009Reactive}\\ 
\citet{Huang2009Dynamic} & Continuum & Pedestrian & Reactive & Yes, 2 & Bi-directional & FDM & Extends the model developed in \citet{Jiang2009Reactive} by introducing look-ahead behavior \\ 
\citet{Jiang2012Numerical} & Continuum & Pedestrian & Reactive & Yes, 2 & Bi-directional & FVM, FSM, TVDRK & Revisits the model developed in \citet{Jiang2009Reactive} and introduces a more general local travel cost\\ 

\hline
\citet{JIANG2011Dynamic} & Continuum & Vehicular & Predictive & No & Isotropic & FVM, FEM, TVDRK & Develops the first PDUE model for vehicular traffic\\
\citet{Du2013Revisiting} & Continuum & Vehicular & Predictive & No & Isotropic & LF, FSM/ SA-MSA & Revisits the route-choice strategy proposed in \citet{JIANG2011Dynamic}\\
\citet{Lin2017Predictive} & Continuum & Vehicular & Predictive & Yes & Isotropic & FVM/ SA-MSA & Extends the model developed in \citet{Du2013Revisiting} to polycentric urban city \\
\citet{Hoogendoorn2004Dynamic} & Continuum & Pedestrian & Predictive & Yes & Isotropic & FDM & Develops the first pedestrian PDUE model comprised of CL and HJB parts\\
\citet{Du2015Reformulating} & Continuum & Vehicular & Predictive & Yes & Anisotropic & LF/ SA-MSA & Revisits the model developed in \citet{Hoogendoorn2004Dynamic} for vehicular traffic and reduces the HJB to HJ\\

\hline
\citet{Xia2009Dynamic} & Continuum & Pedestrian & Predictive \& Reactive & No & Isotropic & DGM, TVDRK& Proposes a route-choice strategy including the global and local route-choice behaviors\\
\citet{Hoogendoorn2015Continuum} & Continuum & Pedestrian & Predictive \& Reactive & Yes & Isotropic & FD & Proposes a more detailed local route-choice formulation compared to \citet{Xia2009Dynamic}\\

\hline
\citet{Jiang2010Higher} & Continuum & Pedestrian & Reactive & No & Isotropic & FVM, FSM, TVDRK & Develops a second-order RDUE model incorporating equilibrium of linear momentum\\
\citet{Jiang2016Macroscopic} & Continuum & Pedestrian & Predictive & No & Isotropic & FVM, FSM/ SA-MSA & Develops a second-order PDUE model similar to \citet{Jiang2010Higher}\\
\citet{Jiang2016Comparison} & Continuum & Pedestrian & Predictive \& Reactive & No & Isotropic & FVM, FSM/ SA-MSA & Compares the second-order RDUE \citep{Jiang2010Higher} and PDUE \citep{Jiang2016Macroscopic} models\\

\hline
\citet{Tao2014Dynamic} & Continuum & Vehicular & --- & No & Isotropic & FVM & Finds the DSO solution instead of the DUE solution\\

\hline
\citet{Laval2017Minimal}  & Continuum and Freeway& Vehicular & ---& Single O-D & Isotropic & Analytical Solution & Investigates the DUE conditions for CS network with an alternative Freeway \\

\hline
\citet{Hanseler2014Macroscopic} & Cells & Pedestrian & Reactive & Yes & Isotropic & Logit & Develops PedCTM model and incorporates cell-based logit route-choice\\
\citet{Hanseler2017Dynamic} & Areas & Pedestrian & Reactive & Yes & Anisotropic & Logit & Replaces the cell-based route-choice in \citet{Hanseler2014Macroscopic} to  link-based route-choice using logit\\

\hline
\citet{Yildirimoglu2014Approximating} & Regions & Vehicular & Predictive & Yes & Isotropic & SNL, MSA, Logit & Divides heterogeneous urban network into homogeneous regions and defines the route as the sequence of regions\\
\citet{Yildirimoglu2015Equilibrium} & Sub-regions & Vehicular & Predictive \& Marginal & Yes & Isotropic & MSA & Divides regions into sub-regions and redefines the route-choice as the sequence of sub-regions\\
\citet{Sossoe2017Reactive} & Zones & Vehicular & Reactive & Yes & Bi-directional & Logit & Decomposes the urban network into zones, which are meshed into cells, and computes the turning proportion in each cell using logit\\
\end{tabular}
\end{center}
\end{table}

\vspace*{-2pt}
\subsection{Departure time choice}

Departure time choice modeling has not been  incorporated in the literature. 
There is a clear need to concentrate in departure time choice within the macroscopic DTA framework, especially after the compelling results by \textbf{\citet{Fosgerau2015Congestion}}, which shows that under "regular sorting" (shorter trips depart later and arrive earlier compared to longer trips) the problem becomes simplified very significantly in the case of a single MFD, to the point that reservoir dynamics do not need to be computed explicitly.
 
 Contrary to the results of  \citet{Fosgerau2015Congestion}, 
\textbf{\citet{Lamotte2017Morning}} show that if the users'
characteristics are such that there is only a single peak in the morning commute, a First-in, First-out (FIFO) sorting pattern emerges within early and late user families having identical $\alpha - \beta - \gamma$ scheduling preferences and heterogeneous trip lengths. The $\alpha,\ \beta$ and $\gamma$ parameters are, respectively, the per-time-unit costs of travel/queuing, earliness, and lateness, as defined first in \citet{vickrey1969congestion}. Empirical measurements are needed to resolve these discrepancies.


\vspace*{-2pt}
\subsection{System optimum}
Only two papers, \citet{Tao2014Dynamic} in continuum space and \citet{Yildirimoglu2015Equilibrium} in discrete space, has tried to establish DSO conditions, but the proposed solution methods are numerical. Since it is hard to extract insights from numerical solution methods, future research should focus on analytical solutions, possibly of simplified systems.

\vspace*{-2pt}
\subsection{Directionality}
Most of the studies propose isotropic models and only a few investigate bi-directional or anisotropic conditions. Although solving anisotropic models seems to be more costly and sometimes unattainable,  it is clear that more effort is needed on investigating anisotropic conditions and the interaction of intersecting flows.

One possibility could be estimating separate MFDs for each cardinal direction, with the method in \citet{Laval2015Stochastic} for example, and incorporate a model to manage the interactions between the flows in different directions.


\subsection{From MFD to the conservation law}

 Surprisingly, we found that among more than 20 studied papers in the continuum-space literature, only three of them have briefly mentioned MFD as a justification for the speed-density relationship \citep{Du2013Revisiting,Du2015Reformulating,Long2017dynamic},  but no attempts can be found in the literature to verify whether or not the assumptions of MFD theory are met. This is
perhaps because in existing models travelers do not exit within the region except at discrete destinations.
But in the case of  origins and destinations distributed randomly in the region, we can postulate the conservation law with source term, i.e.
$\rho_t+\nabla \cdot \mathbf{f}(\rho)=(\lambda-o(\rho))/L$, where $\lambda(x,y,t)$ is the demand flow, which may or may not be constrained by \eqref{capacity constraint}  as discussed earlier\footnote{Notice that in this subsection the density $\rho$ is defined per unit distance, and not per unit area as in the rest of this review.}. In the  isotropic case, we have that the flux vector 
satisfies
$\norm{\mathbf{f}}=V(\rho)\, \rho$, and since by definition (2) we have $o(\rho)=V(\rho)\, \rho\, \frac{L}{\ell}$, it follows that the corresponding conservation law becomes
\begin{equation} \label{Source}
\rho_t+\nabla \cdot \mathbf{f}(\rho)=\lambda/L-V(\rho)\, \rho/\ell. \hskip 4cm \mbox{}
\end{equation}

Notice that  the network length  $L$ and   the trip length $\ell$ appear explicitly in this formulation, unlike existent literature. Research is needed to fully understand  the solution of this equation in the context of the different cost functions proposed in the literature, and to generalize it for the anisotropic case and variable trip length.

\subsection{Capacity constraint }
We saw that a standard assumption in the literature
is   \eqref{capacity constraint}, i.e. that in congestion, the inflow cannot exceed the outflow. 
We argue that this constraint can be relaxed to account for more realistic operations. For example, (i) transient surges in demand could be allowed temporarily, and (ii) distance traveled within the reservoir may increase with congestion as travelers find longer routes.
Research is needed to identify how the capacity constraint should be modified to accommodate these and other modeling improvements.

\subsection{Solution methods }

\subsubsection{Numerical solution of continuum-space models}

For continuum-space models we have seen that the numerical solution methods that have been proposed correspond to standard methods in numerical PDEs. These methods  converge to exact solutions only when the  mesh size of the numerical grid tends to zero;  otherwise,  significant numerical viscosity  can be introduced.
Additionally, it is not clear from the literature how computation times compare with traditional DTA methods; it is possible that for acceptable accuracy  the computational times of continuum-space models might be greater.

Further research is needed to extend the recent advances in \textit{exact }numerical solution methods for one-dimensional kinematic wave (LWR) model \citep{Lighthill1955Kinematic,richards1956shock} to the two-dimensional models in order to develop more efficient and less diffusive numerical methods.
These recent advances exploit (i) a shear symmetry property of conservation laws \citep{Laval2016Symmetries}, and (ii) the link between conservation laws and the Hamilton-Jacobi equation. In the one-dimensional case, the viscosity solution of the Hamilton–Jacobi equation, $\varpi_t + H(\varpi_x)=0$, is
the primitive of the unique entropy solution of the corresponding conservation law, $\rho_t + H(\rho)_x = 0$, where $\rho = \varpi_x$. 

However, in the multi-dimensional case, this one-to-one
correspondence no longer exists, but the gradient $\nabla \varpi$ satisfies a system of conservation laws; see \citet{Kurganov2001Semidiscrete}. Of particular interest would be a cellular automaton type solution method, which could potentially reduce computational times by orders of magnitude.


\subsubsection{Analytical solution of continuum-space models}

As mentioned in section 3.1,  it turns out that the conformal mapping symmetry is valid for all density ranges
and for any MFD.
The implications are profound, as it implies that one can choose the conformal mapping that would simplify boundary conditions as much as possible, thereby simplifying the overall solution method. Research is needed to see if this symmetry is also valid for the other continuum models in the literature.

To show that the conformal mapping symmetry is valid for all density ranges
and for any MFD,  let the discomfort function $g(\rho)=1$ in  Hughes' model without loss of generality, which becomes
\begin{subequations} \label{eqTransform1}
    \begin{empheq}[left = \empheqlbrace\,]{align}
      & \rho_t-\nabla \cdot \left(\rho V^2 \nabla \phi\right)=0, \label{eqTransform1A}\\
       & ||\nabla \phi||=1/V, \label{eqTransform1B}
    \end{empheq}
\end{subequations}
with  appropriate boundary conditions.
We now study \citet{Hughes2002Continuum} equations under a general transformation and then conclude that conformal mappings are the only invariant ones. To see this, let $\mathbf{z}=(x,y)$ represent the original coordinate and let
 \begin{equation} \label{eqTransform3}
	 \mathbf{Z}(x,y)=(X(x,y),Y(x,y)) 
\end{equation}
be a  spatial transformation to the new coordinate system  $\mathbf{Z}=(X,Y)$. The Jacobian 
of this transformation,
\[
\mathbf{J}=\left(
\begin{array}{cc}
 \frac{\partial X}{\partial x} & \frac{\partial X}{\partial y} \\
 \frac{\partial Y}{\partial x} & \frac{\partial Y}{\partial y} \\
\end{array}
\right) 
\]
can be written more concisely as
\begin{equation}
      \mathbf{J}=\nabla_z \mathbf{Z} \label{eqTransform4},
\end{equation}
where we have introduced  the coordinate system as a subscript to  the gradient operators. It can be shown that
\begin{subequations} \label{eqTransform5}
    \begin{align}
& \nabla _z\,\phi (X(x,y),Y(x,y),t)=\nabla _Z\,\phi (X,Y,t) \cdot  \mathbf{J}, \label{eqTransform5A}\\
& \nabla_z\cdot \mathbf{f} (X(x,y),Y(x,y),t)= \nabla_Z\cdot (\mathbf{J}\cdot\mathbf{f}(X,Y,t)). \label{eqTransform5B}
    \end{align}
\end{subequations}

Replacing   \eqref{eqTransform5} in   \eqref{eqTransform1A} gives
\[
\rho_t-\nabla_Z \cdot \left(\rho V^2\mathbf{J}\cdot \nabla_Z \phi\cdot \mathbf{J}\right)=0, 
\]
or,
\[
\rho_t-\nabla_Z \cdot \left(\rho V^2\mathbf{J}\cdot\mathbf{J}^T\cdot\nabla_Z \phi\right)=0. 
\]

Since  $\rho(x,y,t)=\rho(X,Y,t)|\mathbf{J}|$ and $V(\rho(x,y,t))=V(\rho(X,Y,t))|\mathbf{J}|^{-1/2}$ (because transformed small areas equal the original area scaled by the determinant $|J|$), we have
\begin{equation} \label{eqTransform6}
| \mathbf{J}|\rho_t-\nabla_Z \cdot \left(\rho V^2\mathbf{J}\cdot\mathbf{J}^T\cdot\nabla_Z \phi\right)=0.
\end{equation}

Conformal maps have  the unique property that  $ 
\mathbf{J}\cdot  \mathbf{J}^T=| \mathbf{J}|\cdot\mathbf{I}$, 
so that   \eqref{eqTransform6}  becomes
\begin{equation} \label{eqTransform7}
\rho_t-\nabla_Z \cdot \left(\rho V^2\nabla_Z \phi\right)=0.
\end{equation}
which reads exactly as   \eqref{eqTransform1A}. Therefore, conformal mappings are the only ones where Hughes' model is invariant.
%
%

\section*{Acknowledgments} 
This study has received funding from NSF research project \# 1562536.

\section*{References}
\bibliographystyle{elsarticle-harv}
\bibliography{refs}

\end{document}